\documentclass{amsart}
\usepackage{graphicx}
\newtheorem{theorem}{Theorem}

\newtheorem{claim}[theorem]{Claim}

\newtheorem{definition}[theorem]{Definition}

\newtheorem{lemma}[theorem]{Lemma}

\newtheorem{proposition}[theorem]{Proposition}

\newcommand{\abs}[1]{\left\vert#1\right\vert}

\usepackage{amsmath,amscd}
\begin{document}
\title[]{On the Topological Structure Of Complex Tangencies to Embeddings of $S^3$ into $\mathbb{C}^3$}%
\author{Ali M. Elgindi}%

\begin{abstract}
In the mid-1980's, M. Gromov used his machinery of the $h$-principle to prove that there exists totally real embeddings of $S^3$ into $\mathbb{C}^3$. Subsequently, Patrick Ahern and Walter Rudin explicitly
demonstrated such a totally real embedding. In this paper, we consider the generic situation for such     embeddings, namely where complex tangents arise as codimension-2 subspaces. We first consider the Heisenberg group $\mathbb{H}$ and generate some interesting results there-in. Then, by using the biholomorphism of $\mathbb{H}$ with the 3-sphere minus a point, we
demonstrate that every homeomorphism-type of knot in $S^3$ may arise precisely as the set of complex     tangents to an embedding  $S^3 \hookrightarrow \mathbb{C}^3$. We also make note of the (non-generic) situation where complex tangents arise along surfaces.
\end{abstract}

\maketitle


\par\ \par\
\section*{0. Introduction}
In this paper, we will be considering the situation where a 3-manifold is embedded into $\mathbb{C}^3$ in such a way that it assumes its complex tangents along a curve. This is the generic situation for complex tangents to arise in this dimension in a sense that will be described explicitly in the next section. A main result in this subject was proved by Gromov in [4], where he uses the $h$-principle to show that the 3-sphere is the only sphere of dimension bigger than one that admits a totally real embedding into its natural ambient complex space $(S^n \hookrightarrow \mathbb{C}^n)$. Ahern and Rudin subsequently demonstrated an explicit totally real embedding $S^3 \hookrightarrow \mathbb{C}^3$ using a harmonic polynomial in two complex variables (see [1]).
\par\
We use the methods of Ahern and Rudin to consider the generic case where complex tangents arise along knots (and links), and we prove that every topological type of knot may arise as the set of complex tangents to some embedding $S^3 \hookrightarrow \mathbb{C}^3$ (and give notes about possible generalizations to all link types). To arrive at our desired result, however, we will first need to derive analogous results for the Heisenberg group $\mathbb{H} = \{(z, w) | Im(w)= \abs{z}^2\} \subset \mathbb{C}^2$, which is naturally biholomorphic to $S^3$ with the north pole removed.
\par\
The first two sections will focus on preliminaries of complex tangents and knot theory, which then leads to our work on $\mathbb{H}$ in section three, which we then use to achieve our main result for $S^3$ in section four.

\par\ \par\ \par\ \par\ \par\
\section*{1. Preliminaries of Complex Tangents}
Let $M$ be a real manifold of dimension $k$ and suppose we embed $M \subset \mathbb{C}^n$ in a smooth manner (in fact, we need only a $\mathcal{C}^1$-embedding). We say a point $x \in M$ is $\emph{complex tangent}$ if the tangent space to $M$ at $x$ contains a complex linear subspace. At such a point $x$, we will have
$T_x (M) \cap J(T_x (M)) \neq \{0\}$, where $J: T_* (\mathbb{C}^n) \rightarrow T_* (\mathbb{C}^n)$ is the isomorphism given by multiplication with $i$ (on each tangent space).
\par\
Let $\aleph_M$ be the set of complex tangents to $M \subset \mathbb{C}^n$. Note that at each complex tangent point $x \in \aleph_M$, the tangent space $T_x (M)$ contains a (non-trivial) linear complex subspace, and as such has a "maximal" complex subspace, namely: $T_x (M) \cap J(T_x (M))$. The (complex) dimension of this subspace is called the $degree$ of the complex tangent $x$.
\par\
The topological structure of $\aleph_M$ in general dimensions can be immensely complicated, with degrees varying amongst the tangents; the global structure can be very singular and stratified.
\par\ \par\
Of much interest is the question of existence of $\emph{totally real}$  embeddings into $\mathbb{C}^n$, i.e. the situation where $\aleph_M = \emptyset$. By mere dimensionality observations, we see this situation is possible only when $k \leq n$, i.e. the real dimension of the manifold is less than or equal the complex dimension of the ambient space. Note that $k=n$ is the topologically natural situation, as every real $n$-manifold may be embedded in $\mathbb{R}^{2n}$ (=$\mathbb{C}^n$) by the Whitney embedding theorem. It is this situation that is of interest to us.
\par\
We note that not every manifold may be embedded in a totally real manner, in particular no sphere of any dimension except one or three may be embedded totally real. This was proved by M. Gromov (see [4]). In dimension one this result is trivial, while in dimension three Ahern and Rudin in [1] demonstrate for us an explicit totally real embedding of $S^3$.
\par\ \par\
In this paper we focus on the case where $n=3$, i.e. $M^3 \hookrightarrow \mathbb{C}^3$. Note that here we consider embeddings as maps of an intrinsic manifold. In this situation,  there is only one class of complex tangents, that is, they must all have degree 1. This is clear as $T_x (M)$ is a real 3-dimensional vector space, and hence can contain at most a complex line. Now, consider the set $\textbf{\emph{G}}_{6,3} =\{3$-planes $P \subset \mathbb{R}^{6} \}$. Note $\mathbb{R}^{6} \cong \mathbb{C}^{3}$, so we have the isomorphism $J$ on the vector space.
\par\
We also consider the subset consisting of "partially complex" 3-planes, which we call:
\par\ \par\
$\mathbb{W} = \{P \in \textbf{\emph{G}}_{6,3} | P \cap J(P) \neq \emptyset\} \subset \textbf{\emph{G}}_{6,3}$
\begin{lemma} $\mathbb{W}\subset \textbf{\emph{G}}_{6,3}$ is a smooth submanifold of codimension 2. \end{lemma}
\par\
$\textbf{\emph{\underline{Proof}:-}}$
\par\
Let's demonstrate the dimensionality relation in general dimension $n$. Let $\mathbb{W}_1 = \{P \in \textbf{\emph{G}}_{2n,n} | P \cap J(P) \cong \mathbb{C}\}$ be the set of all $n$-planes in $\mathbb{C}^n$ containing exactly a complex line as its maximal complex subspace. Choose a specific complex line $L \in \mathbb{C}\textbf{\emph{G}}_{n,1}$, or equivalently, $L \in \textbf{\emph{G}}_{2n,2}$ so that $L \bigcap J(L) = L$. Note that a generic $(n-2)$-plane $Q \subset  L^\perp \cong \mathbb{R}^{2n-2}$ will be totally real. As a result, for almost any such $(n-2)$-plane $Q$, the $n$-plane $Q \bigoplus L \in \textbf{\emph{G}}_{2n,2}$ will satisfy that $Q \bigoplus L \bigcap J(Q \bigoplus L) = L$.
\par\ \par\
Recall that $\mathbb{C}\textbf{\emph{G}}_{n,1} = \mathbb{C}\mathbb{P}^{n-1}$, which is a complex manifold of complex dimension $n-1$. By our note above, the dimension of all $(n-2)$-planes $Q$ which are totally real is $dim(\textbf{\emph{G}}_{2n-2,n-2}) = (n-2) (2n-2 -(n-2)) = n(n-2)$. As every element of $\mathbb{W}_1$ may be uniquely expressed as $Q \bigoplus L$ for a complex line $L$ and such a totally real $(n-2)$-plane $Q$, we find that dimension of $\mathbb{W}_1$ must be
 \par\
$dim_{\mathbb{R}}(\mathbb{W}_1) = n(n-2) + 2n-2 = n^2 -2$.
\par\
As $dim(\textbf{\emph{G}}_{2n,n}) = n^2$, we get the desired result that $\mathbb{W}_1 \subset \textbf{\emph{G}}_{2n,n}$ is of codimension 2.
\par\ \par\
Furthermore, for $n = 3$, we see that $\mathbb{W} = \mathbb{W}_1 \subset \textbf{\emph{G}}_{6,3}$ is a smooth submanifold as it arises precisely as an orbit to the natural group action of the compact group $U(3)$ on $\textbf{\emph{G}}_{6,3}$.
\par\ \par\
$\textbf{\emph{QED}}$
\par\ \par\
For manifolds of general dimension, the space $\mathbb{W}$ may be very singular, and we will not investigate this general case.
\par\ \par\
 In the following theorem, we will demonstrate by using the Trasversality Theorem that any (smooth) embedding of a given 3-manifold $M \subset \mathbb{C}^{3}$  may be perturbed an arbitrarily small amount so that it will have its set of complex tangents (if any), $\aleph_M \subset M$  arising as a submanifold of codimension 2; in other words, along a curve. This result is in fact a special case of a theorem of Webster, which is stated for general dimensions in his paper [10].
 \par\
 We will denote by $G:M \rightarrow \textbf{\emph{G}}_{6,3} $ the Gauss map of an embedding $M \subset \mathbb{C}^3 \cong \mathbb{R}^6$, given by: $G(x) = T_x (M) \subset \mathbb{R}^6$, for any $x \in M$. Then by definition, a point $x \in M$ is complex tangent if and only if its image under the Gauss map is contained in $\mathbb{W}$.
\par\
\begin{lemma} For any embedded $M \subset \mathbb{C}^3$, there exists an arbitrarily small perturbation of the embedding whose Gauss map is transverse to $\mathbb{W}$ \end{lemma}
\par\ \par\
$\textbf{\emph{\underline{Proof}:-}}$
\par\
Consider the map $F: M \times SO(6)  \rightarrow \textbf{\emph{G}}_{6,3}$ which is defined by:
\par\
$F(x,l) = F_l (l \cdot x)$, where $\{F_l : l \cdot M  \rightarrow \textbf{\emph{G}}_{6,3}\}$ is the collection of Gauss maps for the embeddings given by the rotations of  $M \subset \mathbb{C}^3$ by elements $l \in SO(6)$. By linearity, it is easy to see that $F(x, l) = l \cdot T_x (M)$, for $x \in M$ and $l \in SO(6)$. We note that  for any given plane $F(x)$, by rotating the plane by the elements of $SO(6)$ we can obtain any element of the Grassmannian $\textbf{\emph{G}}_{6,3}$. Now consider the differential of $F$:
\par\
$dF : T_x (M) \times T_l (SO(6))  \rightarrow  T_{l \cdot T_x(M)} (\textbf{\emph{G}}_{6,3})$. Since every element of the Grassmannian is assumed by rotating $F(x)$ by the elements $l \in SO(6)$ (for any $x \in M$), we see that every tangential direction in the Grassmannian (at any given point) can be taken as the derivative of some path in $M \times SO(6)$, in particular every element of $T_{l \cdot T_x(M)} (\textbf{\emph{G}}_{6,3})$ is assumed by $dF$. From this, we see that $dF$ must be onto.
\par\ \par\
As a result, the map $F$ is transverse to the subset $\mathbb{W} \subset \mathbb{C}^3$; this is trivial as $Image(dF) = T (\textbf{\emph{G}}_{6,3})$ at every point. Now, by applying the Parametric Transversality Theorem (see [5]), we see that for any embedded 3-manifold $M$, the embedding $l \cdot M$ given by rotating the given embedding by a matrix $l \in SO(6)$ will be transverse to $\mathbb{W}$, for almost any element $l \in SO(6)$. In particular, given any embedding there exists an arbitrarily small perturbation, given by rotating with a matrix arbitrary close to the identity, so that the perturbed embedding is transverse to $\mathbb{W}$.
\par\ \par\
$\textbf{\emph{QED}}$
\par\ \par\
Using the normal version of the Transversality Theorem (see [5]) and with the above lemma, we arrive at the following result:
\par\
\begin{proposition} For any embedded $M \subset \mathbb{C}^3$, there exists an arbitrarily small perturbation of the embedding so that the embedding is either totally real or takes its complex tangents along a smooth curve (or curves).   \end{proposition}
\par\ \par\
$\textbf{\emph{\underline{Proof}:-}}$
\par\
We proved earlier in the section that $\mathbb{W} \subset \textbf{\emph{G}}_{6,3}$ is codimension two. Furthermore, by the above lemma there exists an arbitrarily small perturbation so that the Gauss map of the embedding is transverse to $\mathbb{W}$. Hence, we may conclude by the Transversality Theorem that if the image of the Gauss map intersects $\mathbb{W}$, the inverse image of $\mathbb{W}$ must also be of codimension 2 in $M$. As $\aleph_M = G^{-1} (\mathbb{W})$ is the set of complex tangents and $M$ is 3-dimensional, the proposition follows.
\par\ \par\
$\textbf{\emph{QED}}$
\par\ \par\
We can now consider embeddings of 3-manifolds whose complex tangents arise along curves to be $\emph{generic}$. In particular, if $M$ is closed, a generic embedding of $M$ will assume its complex tangents along knots, or more generally along links of knots. This terminology of being generic is used, for example, by Webster in [10].
\par\ \par\
We will now restrict our consideration to the subcollection of 3-manifolds that arise as smooth submanifolds of $\mathbb{C}^2$, i.e. $M = \{\rho (z, w) = 0\}$, where $\rho:\mathbb{C}^2 \rightarrow \mathbb{R}$ is a smooth map whose differential is never zero on $M$. For such manifolds, embeddings into $\mathbb{C}^3$ arise very
naturally as graphs of maps $f: \mathbb{C}^2 \rightarrow \mathbb{C}$, which are sufficiently smooth on a neighborhood of $M$.
\par\
In particular, let $f:\mathbb{C}^2 \rightarrow \mathbb{C}$ be such a map and let $F:M \hookrightarrow  \mathbb{C}^3$ be its graph, $F=graph(f|_M )$. By definition, $F(z, w) = (z, w, f(z, w))$ and is an embedding by virtue of construction (as a graph).
\par\ \par\
On such a manifold $M = \rho^{-1} (0)$, there exists a operator which we call the $\emph{CR-operator}$, and is given by:
\par\
$\mathcal{L}_M = \frac{\partial\rho}{\partial\overline{w}}  \frac{\partial}{\partial\overline{z}} -
\frac{\partial\rho}{\partial\overline{z}}  \frac{\partial}{ \partial \overline{w}}$,  acting on complex-valued functions defined in a neighborhood of $M$.
\par\ \par\
We will now outline the work of Ahern and Rudin which they published in [1]. They use the following proposition:
\par\
\begin{proposition} Let $M = \rho^{-1} (0) \subset \mathbb{C}^2$ be a real hypersurface as above. Given a graphical embedding $F(z, w) = (z, w, f(z,w)): M \hookrightarrow \mathbb{C}^3$, the embedding will have its set of complex tangents precisely where the tangential Cauchy-Riemann operator applied to f is zero (on $M$). \end{proposition}
\par\ \par\
The statement in the above proposition may be written more precisely:
\par\ \par\
$\aleph_f = \aleph_{F(M)} = \{(z,w) \in \mathbb{C}^2 | \mathcal{L}_M (f) (z, w) = 0$ and $\rho (z,w) = 0\} \subset M$
\par\ \par\
In the case of $M=S^3$, we have $\rho(z,w) = \abs{z}^2+\abs{w}^2-1$ and so: $\mathcal{L}_{S^3} = \mathcal{L} = w \frac{\partial}{\partial\overline{z}}-z\frac{\partial}{\partial\overline{w}}$.
\par\
In general, this operator (over $S^3$) is very difficult to analyze, and has been the subject of research
for many years by many great mathematicians. The operator on the Heisenberg group turns out to be a little simpler and much more manageable.
\par\ \par\
Note that the potential configurations for $\aleph_f$ are numerous, in particular  $\aleph_f$ could be empty, discrete points, curves, surfaces, or unions of such.
Note also if $f$ is holomorphic, $\aleph_f = M$ (every point is complex tangent).
\par\
As we saw above, the generic situation will be that $\aleph_f$ is some curve, or union of curves. As we will demonstrate, all different kinds of topological configurations are possible, both in the generic situation and the non-generic situation. In particular we show this
for the Heisenberg group $\mathbb{H}$ and the three-sphere $S^3$ (both submanifolds of $\mathbb{C}^2$).
\par\
Before we can proceed further, we will need to address some relevant questions regarding the topology of knots in $\mathbb{R}^3$ and $S^3$.
\par\ \par\

\section*{2. Relevant Notes on Knot Topology}
A $\emph{knot}$ is a smooth (at least continuous) simple embedding of the circle $S^1$ (no double points). In $\mathbb{R}^2$ they are closed Jordan curves. Note that the figure 8 curve (in $\mathbb{R}^2$) is not a knot as it has a double point. Of great interest is the structure of such curves in $\mathbb{R}^3$. However, both practically and traditionally, the natural ambient space to consider knots is $S^3$. We note that embedding a closed curve in a higher dimensional space (such as $\mathbb{R}^4$) is trivial as there is enough dimensions to "unravel any knot", and so all such embeddings are topologically circles.
\par\
A $\emph{link}$ is defined as a disjoint union of knots... note there are further invariants for links given by $\emph{linking numbers}$.
\par\
For a complete exposition of the theory of knots, we refer the reader to Lickorish ([6]).
\par\ \par\
In $S^3$, we may classify knots up to diffeomorphism (homeomorphism):
\begin{definition} We say two smooth (continuous) knots  $K_1, K_2 \subset S^3$ are
$\emph{topologically}$ $\emph{equivalent}$, or of the same $\emph{topological type}$, if there exists a diffeomorphism (homeomorphism) $h:S^3 \rightarrow S^3$ so that $h(K_1)=K_2$.    \end{definition}
\par\ \par\
This gives us an equivalence relation on the set of knots, and note that the pair: $(S^3, S^3 \setminus K)$ is now a topological invariant on knots (links) $K \subset S^3$. Let $\lambda$ be an arbitrary equivalence class of knots (in $S^3$) under this relation. We say a knot $K \subset S^3$ is of $\emph{topological type}$ $\lambda$ if $K$ is an element of the equivalence class $\lambda$.
\par\
Note we may extend this definition of topological type to general links, although we will have to keep track of the linking numbers of the individual knot components.
\par\ \par\
We can classify (most) knots in $\mathbb{R}^3$ directly from the classification of knots in $S^3$; in particular let $\sigma:S^3 \setminus \{pt\} \rightarrow \mathbb{R}^3$ be a stereographic projection through a given point in $S^3$, and let $\kappa$ be its inverse. There are two inherent types of knots in $ \mathbb{R}^3$: first, we say a knot in $\mathbb{R}^3$ is $\emph{bounded}$ if it is a simple closed curve. We say a knot in $\mathbb{R}^3$ is $\emph{unbounded}$ if it is a curve that does geometrically knot, but has "two rays to infinity"; one could see Shastri work in [9] for further description of such knots.
\par\ \par\
We may then classify bounded knots in $\mathbb{R}^3$ in the following manner:
\par\
\begin{definition} A bounded knot $K \subset \mathbb{R}^3$ is said to be of $\emph{bounded type}$ $\lambda$ if its image: $\kappa(K) \subset S^3$ is of topological type  $\lambda$ as a knot in $S^3$. \end{definition}
\par\ \par\
Hence, bounded knots in $\mathbb{R}^3$ are classified precisely by the classification of their image under stereographic projection diffeomorphism; its type exactly corresponds to the topological type of the corresponding knot in $S^3$. One easily sees that this extends to an equivalence relation and topological invariant of bounded knots in $\mathbb{R}^3$.
\par\
We may classify (certain) unbounded knots in $\mathbb{R}^3$ in a similar way:
 \begin{definition} Let $K \subset \mathbb{R}^3$ be an unbounded knot. We say $K$ is of $\emph{unbounded type}$ $\lambda$ if the image of K under the stereographic projection, $\kappa(K)$ forms a continuous (in particular, well-defined) knot if we were to add the "point at infinity", and this union assumes the topological type $\lambda$ as a knot in $S^3$. \end{definition}
 \par\ \par\
 We see again that these unbounded knots in $\mathbb{R}^3$ are classified precisely by the topological type of their images in $S^3$.
\par\
Note that only certain unbounded knots $K \subset \mathbb{R}^3$ can arise in this fashion, i.e. as the image of a knot in $S^3$ under the projection $\sigma$. We define such knots to be unbounded knots of $\emph{finite count}$. More precisely, we say an unbounded knot in $\mathbb{R}^3$ is of finite count if there exists a positive number $R>0$ so that outside the ball of radius $R$ about the origin the knot consist precisely of two "unknotted" rays going to infinity. More precisely, if the complement: $K \bigcap (\mathbb{R}^3 \setminus \textbf{B}_R (\overrightarrow{0}))$ is topologically equivalent to the set: $\{(0, 0, t) | t \geq 1\} \bigcup \{(0, 0, t) | t \leq -1\}$. Equivalently, we can define an unbounded knot as being of finite count if its projections onto the coordinate planes has finitely many double points.
\par\
If the unbounded knot is not of finite count (i.e. there exists no such number $R$, or equivalently its projections admit infinitely many crossings), we say the knot is of $\emph{infinite count}$. Note that the image of such an unbounded knot $K$ of infinite count under stereographic projection, $\kappa(K)$, cannot converge at the point $(0,1)$, and that the image of any knot on $S^3$ under the inverse stereographic projection must be of finite count (by properties of knots as submanifolds). Hence, for our purposes we will only interested in the unbounded knots that are of finite count. Therefore, we will essentially ignore unbounded knots of infinite count for the remainder of this paper, and make no claims there-in regarding such knots.
\par\ \par\
We may further generalize the above notions to general links. The topological type (or classification) of a link will depend only on the topological type of its knot components and the linking numbers between the knots there-in.
In most of what follows, we will limit our consideration to only single knots, however we may be able to directly generalize some results to $\emph{simple}$ links, that is links whose knot components have linking numbers at most one with each other.
\par\ \par\
Although the construction above suits well enough for an investigation into the topological classification of knots (on $S^3$ and so by extension on $\mathbb{R}^3$), there are many perspectives and tools for investigating knots. One such perspective is old and direct, where we classify knots by their projections onto the coordinate planes. While this perspective is interesting and intuitive, it serves to be tedious and unhelpful for our purposes. There are also numerous topological invariants, such as the Jones and Alexander polynomials, which are helpful in understanding knots, but are also immensely complicated and not directly applicable in our case (as far we know).
\par\
We will not be making use of heavy machinery or such perspectives in this paper, and in fact wish to simplify our computations as much as possible. As we mentioned before, we wish to show that every topological type of knot may occur precisely as the set of complex tangents to an embedding, and we will do so through the means of complex polynomials and algebraic sets. The result of S. Akbulut and King (see [2]) gives us that "all knots are algebraic" in $S^3$, or more precisely:
\begin{theorem}: (Akbulut-King) Let $K \subset S^3 \subset \mathbb{R}^4$ be a knot. Then there exists an algebraic set $\mathfrak{P} \subset \mathbb{R}^4$ (in two real equations) passing through the origin and only singular at the origin, so that the knot $\widetilde{K} = \mathfrak{P} \bigcap S^3$ is of the same topological type as $K$, i.e. $(S^3, S^3 \setminus K) \cong (S^3, S^3 \setminus \widetilde{K})$ (as topological pairs).   \end{theorem}
\par\ \par\
(This is in fact a weaker version of the theorem proved by Akbulut and King... please see reference for the more general results).
\par\ \par\
In particular, this theorem shows us that every topological type (or class) of knots in $S^3$ has an algebraic representative. In other words, given any knot $K$, we may assume (up to topological type):
\par\
$K = \{p=0, q=0\} \bigcap S^3$, where $p,q:\mathbb{R}^4 \rightarrow \mathbb{R}$ are real polynomials. In particular, we may take:
$f= p+iq:\mathbb{C}^2 \rightarrow \mathbb{C}$ and have: $K = \{p=0, q=0\} \bigcap S^3 = \{f=0\} \bigcap S^3 \subset \mathbb{C}^2 = \mathbb{R}^4$.
\par\
Hence, we may rephrase the above theorem as follows: for any knot $K \subset S^3$, there exists a complex polynomial (in $z, w, \overline{z}, \overline{w}$ coordinates), call it $f$, so that:
$K \cong \widetilde{K} = \{f=0\} \bigcap S^3$.
\par\
Further, $\widetilde{K} \subset S^3$ will be non-singular as a real algebraic variety.
\par\ \par\
$\textbf{\emph{\underline{Remark:}}}$ In fact, as any link may also arise as the boundary of a Seifert surface, the proof of Akbulut-King in their paper seems to generalize to any equivalence class of a general link. See the reference [2].
\par\ \par\
We will generalize (or rather, restrict) the result of Akbulut-King (Theorem 8) to $\mathbb{R}^3$. However, we will first investigate the Heisenberg group $\mathbb{H}$ (which is diffeomorphic to $\mathbb{R}^3$) to avoid confusing notation and allow for a more natural construction. We shall in fact consider $\mathbb{H}$ and $\mathbb{R}^3$ as interchangeable via the standard diffeomorphism (see next section).
\par\ \par\
We will also note, by [8], that we can similarly give every topological type of surface in $\mathbb{R}^3$ (or $S^3$) a polynomial representative.
\par\ \par\
\section*{3. Analysis on the Heisenberg Group}
In this section, we will define the Heisenberg group and consider knots in this space. The analysis of complex tangents to embeddings of this space will be easier to study and we will be able to obtain some interesting results there-in.
\par\ \par\
Consider the complex space $\mathbb{C}^2$ with holomorphic coordinates $z, w$. The Heisenberg group is a (real) hypersurface of $\mathbb{C}^2$ given by: $\mathbb{H} = \{(z, w) | Im(w)= \abs{z}^2\}$. Note that $\mathbb{H} \cong \mathbb{R}^3 \cong \mathbb{C} \times \mathbb{R}$ are differomorphic via:
$(z, w) \rightarrow (z, u)$, where $u=Re(w)$.
\par\ \par\
In complex coordinates, one easily finds that the tangential Cauchy-Riemann operator to $\mathbb{H}$ is given by:
\par\
$\mathcal{L}_\mathbb{H} = \frac{\partial\rho}{\partial\overline{w}}  \frac{\partial}{\partial\overline{z}}
- \frac{\partial\rho}{\partial\overline{z}}  \frac{\partial}{ \partial \overline{w}} = 2z  \frac{\partial}{\partial\overline{w}} + i \frac{\partial}{\partial\overline{z}}$, as here $\rho(z,w)=i(\overline{w}-w)-2z\overline{z}$
\par\ \par\
By our earlier discussion, for any (sufficiently) smooth map $f: \mathbb{C}^2 \rightarrow \mathbb{C}$, the zeros of the function:
$\mathcal{L}_\mathbb{H}(f)|_\mathbb{H}$ give precisely the set of complex tangents to the embedding: $F:\mathbb{H} \hookrightarrow \mathbb{C}^3$, where $F=graph(f|_\mathbb{H})$.
\par\ \par\
We immediately find many (even linear) totally real embeddings of $\mathbb{H}$, in particular the graph of $f(z,w)=\overline{z}$ would be totally real, with  $\mathcal{L}_\mathbb{H} (f) = i$, a (non-zero) complex constant.
\par\ \par\
The question of the solvability of this differential equation, namely: $\mathcal{L}_\mathbb{H} (f) = h$, for a given $h:\mathbb{C}^2 \rightarrow \mathbb{C}$, is generally very difficult even in the local sense. In fact, Hans Lewy demonstrated a complex function for which there exists no solution.
\par\
We find however, to our pleasant surprise, that this question is readily solvable for polynomials (in the global sense), as exhibited in the following lemma:
\begin{lemma} As a linear operator on infinite dimensional space of all complex polynomials, $\mathcal{P}$,  $\mathcal{L}_\mathbb{H}$ is onto; i.e.:  $\mathcal{L}_\mathbb{H}: \mathcal{P} \rightarrow \mathcal{P}$ is a surjective linear map.            \end{lemma}
\par\ \par\
$\textbf{\emph{\underline{Proof}:-}}$
 \par\
 We first note that by the structure of $\mathcal{L}_\mathbb{H}$ as a linear differential operator with polynomial (linear and constant) coefficients, that $\mathcal{L}_\mathbb{H}$  indeed maps polynomials to polynomials. In fact,  $\mathcal{L}_\mathbb{H}$ will map polynomials of (homogeneous) weight $s$ to polynomials of weight $(s-1)$, where $z, \overline{z}$ have weight 1, and $w$ (or $u$) has weight 2.
\par\
Our goal is to show that every polynomial may arise in the image of the operator. As  $\mathcal{L}_\mathbb{H}$ is linear, it will suffice to show all monomials are in the range. Let $f(z,w)=z^j \overline{z}^k w^m \overline{w}^l$ be such a monomial... note that as $\mathcal{L}_\mathbb{H} (w^m g(z,w))= w^m \mathcal{L}_\mathbb{H}(g(z,w))$, we need only verify the identity for all monomials of the form: $f(z,w)=z^j \overline{z}^k \overline{w}^l \in \mathcal{L}_\mathbb{H}(\mathcal{P})$.
\par\
We proceed by induction on $l \in \mathbb{N}$. Note if $l=0$, then $f(z,w)=z^j \overline{z}^k = \mathcal{L}_\mathbb{H}
(\frac{1}{i(k+1)} z^j \overline{z}^{k+1})$, for any $j, k \in \mathbb{N}$.
\par\
Now, let $l \in \mathbb{N}$ be fixed and suppose: $z^j \overline{z}^k \overline{w}^r \in \mathcal{L}_\mathbb{H}(\mathcal{P})$, for all $j, k \in \mathbb{N}$, and $r\leq l$.
\par\
We wish to show $f(z,w)=z^j \overline{z}^k \overline{w}^{l+1} \in \mathcal{L}_\mathbb{H}(\mathcal{P})$. We have that
\par\
$\mathcal{L}_\mathbb{H}(\frac{1}{i(k+1)} z^j \overline{z}^{k+1} \overline{w}^{l+1}) = z^j \overline{z}^k \overline{w}^{l+1} + \frac{2(l+1)}{i(k+1)} z^{j+1} \overline{z}^{k+1} \overline{w}^l$.
\par\
Note that the second term on the right-hand side is in the range of $\mathcal{L}_\mathbb{H}$ by the induction hypothesis and linearity... Let $g: \mathbb{C}^2 \rightarrow \mathbb{C}$ be a polynomial such that $\mathcal{L}_\mathbb{H}(g)= \frac{2(l+1)}{i(k+1)} z^{j+1} \overline{z}^{k+1} \overline{w}^l$
The term on the left hand side is in the range by construction. Hence, we find that
\par\
$f(z,w)=z^j \overline{z}^k \overline{w}^{l+1}= \mathcal{L}_\mathbb{H}(\frac{1}{i(k+1)} z^j \overline{z}^{k+1} \overline{w}^{l+1})
 - \mathcal{L}_\mathbb{H}(g) \in \mathcal{L}_\mathbb{H}(\mathcal{P})$
 \par\
 by the linearity of the operator.
\par\
Therefore, by the principle of induction and our preliminary arguments, we have shown that $\mathcal{L}_\mathbb{H}(\mathcal{P}) = \mathcal{P}$, and our lemma is proven.
\par\ \par\
$\textbf{\emph{QED}}$
\par\ \par\
We have shown that every polynomial is in the range of the CR-operator (over $\mathbb{H}$). In particular, given any algebraic set in $\mathbb{H}$ which is the zero set of two real polynomial equations, there exists an embedding of $\mathbb{H} \hookrightarrow \mathbb{C}^3$ whose complex tangents are precisely that algebraic set.
\par\
We deduce this more formally: given any such algebraic set on $\mathbb{H}$, we can write it as the intersection of the zero set of a complex polynomial $p$ with $\mathbb{H}$. Here $p = f+ ig$ and $f=0, g=0$ define the algebraic set.
\par\ \par\
As $p \in \mathcal{P} = \mathcal{L}_\mathbb{H}(\mathcal{P})$, there exists a complex polynomial $q$ so that $\mathcal{L}_\mathbb{H}(q) = p$. Then the set of complex tangencies to the map $F:\mathbb{H} \hookrightarrow \mathbb{C}^3$, given by $F=graph(q|_\mathbb{H})$, is precisely
\par\
$\aleph_q = \{\mathcal{L}_\mathbb{H}(q) (z,w) = 0\} \bigcap \mathbb{H} = \{p(z,w)= 0\} \bigcap \mathbb{H}$
 \par\
 The set $\aleph_q$ is precisely our given algebraic set!
\par\ \par\
Now, we wish to demonstrate that (almost) every topological class of knots in $\mathbb{H} \cong \mathbb{R}^3$ admits an algebraic representative. We will prove this by "restricting" the (known) corresponding result for $S^3$ given by Akbulut-King.
Note we will need to show this for both unbounded and bounded knot types.
\par\ \par\
Recall that in the previous section we classified all knots in $\mathbb{R}^3$ which are either bounded or unbounded of finite count using the classification of knots on $S^3$. As we noted before, $\mathbb{R}^3 \cong \mathbb{H}$ are naturally diffeomorphic, and hence we may classify all such knots in $\mathbb{H}$ in an analogous manner.
\par\
In particular, let $\varphi: \mathbb{H} \rightarrow S^3 \setminus \{(0,1)\}$ be the standard biholomorphism; namely:
$\varphi(z,w) = (\frac{2z}{w+i},\frac{w-i}{w+i})$, which is smooth when restricted to $\mathbb{H}$.
\par\
Let $\psi: S^3 \setminus \{(0,1)\} \rightarrow \mathbb{H}$ be its inverse, in particular:
$ \psi(z,w) = (\frac{iz}{1-w},i\frac{1+w}{1-w})$.
\par\ \par\
Note that $\psi$ can be considered as a natural stereographic projection (again considering $\mathbb{H}$ and $\mathbb{R}^3$ as interchangeable spaces), and we will use this $\psi$ to formulate the classification of knots (and links) on $\mathbb{H} \cong \mathbb{R}^3$ in the manner we demonstrated in the previous section. We will continue this particular formulation (in terms of $\psi, \varphi$) for the rest of the paper.
\par\ \par\
In other words, for a bounded knot  in $\mathbb{H}$ its topological type is determined completely by its image under $\varphi$ and similarly we may classify unbounded knots of finite count, which are now defined analogous to our definition before (on $\mathbb{R}^3$). We will now use this formulation and make it more precise in the proof of the following theorem:
 \par\
\begin{theorem} (All knots are algebraic in $\mathbb{H}$)
Let $K \subset  \mathbb{H}$ be a knot of bounded type or unbounded of finite count. Then there exists a real algebraic set $\mathfrak{P} \subset \mathbb{R}^4 \cong \mathbb{C}^2$ (in two real equations) so that the knot $\widetilde{K} = \mathfrak{P} \bigcap \mathbb{H}$ is of the same topological type as $K$.
\end{theorem}
\par\ \par\
$\textbf{\emph{\underline{Proof}:-}}$
\par\
Let $K \subset \mathbb{H}$ be a knot of bounded type, that is $K$ is a simple closed curve. Then $\varphi(K) \subset S^3$ is also a knot (of corresponding smoothness) and hence has a topological type, call it $\lambda$. Then we can say $K \subset \mathbb{H}$ is of bounded topological type $\lambda$. Now, by the Theorem of
Akbulut and King, there exists a knot $\widetilde{K} \subset S^3$ that is algebraic, say $\widetilde{K} = \{p=0\} \bigcap S^3$, for some complex polynomial $p$, so that $\widetilde{K} \cong \varphi(K)$. That is that $\widetilde{K}$ and $\varphi(K)$ have same topological type (in $S^3$), which was called $\lambda$.
\par\ \par\
We may assume without loss of generality that $\widetilde{K} \subset S^3$ does not pass through (0,1), by rotating the sphere (linearly) to ensure this.
Consider then $\psi(\widetilde{K}) \subset \mathbb{H}$. As $\widetilde{K}$ misses the point at infinity, $\psi(\widetilde{K})$ will be a bounded knot of type $\lambda$ (as $\psi$ is a biholomorphism).
\begin{claim} $\psi(\widetilde{K}) \subset \mathbb{H}$ is algebraic.  \end{claim}
\par\
It is quite clear that birational maps preserve algebraic sets, but let's demonstrate this.
We know that $\widetilde{K} = \{p=0\} \bigcap (S^3 \setminus \{(0,1)\})$, where $p: \mathbb{C}^2 \rightarrow \mathbb{C}$ is a complex polynomial.  Hence, $\psi(\widetilde{K}) = \psi(\{p=0\} \bigcap S^3 \setminus \{(0,1)\}) =\psi(\{p=0\}) \bigcap \psi(S^3 \setminus \{(0,1)\}) = \psi(p^{-1}(0)) \bigcap \mathbb{H}$, as $\psi$ is a biholomorphism.
\par\ \par\
We claim that $\psi(p^{-1}(0)) = q^{-1}(0)$, where $q$ is the numerator of the (rational) function: $p \circ \varphi:\mathbb{H} \rightarrow \mathbb{C}$. A simple computation gives
\par\
$p(\varphi(z,w))=\frac{q(z,w)}{(i+w)^M (\overline{w}-i)^N}$
\par\
for some $M, N \in \mathbb{N}$ and $q$ a complex polynomial. Note that $q$ is unique up to multiplicative factors of $(i+w)$ and its conjugate, which are never zero on $\mathbb{H}$. Hence, we have that $(p \circ \varphi)^{-1} (0) = q^{-1} (0)$.
\par\ \par\
But then: $ q^{-1} (0) \bigcap \mathbb{H}= (p \circ \varphi)^{-1} (0) \bigcap \mathbb{H}=
\varphi^{-1} (p^{-1} (0))\bigcap \mathbb{H}$
\par\
$= \psi (p^{-1} (0))\bigcap \mathbb{H} = \psi (p^{-1} (0) \bigcap S^3 \setminus \{(0.1)\})  = \psi(\widetilde{K})$, and thus our claim is proven.
\par\ \par\
Hence, $\overbrace{K}=\psi(\widetilde{K}) \subset \mathbb{H}$ is an algebraic bounded knot topologically equivalent to the original (arbitrary) knot $K \subset \mathbb{H}$ (following the biholomorphisms). Hence, every bounded knot type in $\mathbb{H}$ has an algebraic representative. We also get the analogous result for $\mathbb{R}^3$ (as $\mathbb{R}^3 \cong \mathbb{H}$ is polynomial).
\par\ \par\
We get the analogous result for unbounded knots of finite count: let $K \subset \mathbb{H}$ be any unbounded knot of finite count. Then $\varphi(K) \cup \{(0,1)\} \subset S^3$ will be a (well-defined) knot passing through $(0,1)$; suppose it is of type $\lambda$ as a knot in $S^3$.
Then, again by Akbulut-King, there exists a algebraic knot of type $\lambda$ (i.e., equivalent to $\varphi(K)$), call it $\widetilde{K} \subset S^3$. By the use of a linear rotation, we may assume without loss of generality that $\widetilde{K}$ contains $(0,1)$. By our above arguments, $\psi(\widetilde{K} \setminus \{(0,1)\}) \subset \mathbb{H}$ will be an algebraic (unbounded) knot of the same type as $K$, and as $K$ was chosen arbitrarily, we again have that every topological type of unbounded knot (of finite count) in $\mathbb{H}$ admits an algebraic representative. We also get the analogous result for such knots in $\mathbb{R}^3$.
\par\ \par\
Hence, we find that every knot type (which is bounded or unbounded of finite count) may arise as an algebraic subset of $\mathbb{H}$. In particular, choose any equivalence class of knots $\lambda$ (bounded or unbounded). Then there exists some complex polynomial $p$ so that:
$K = p^{-1} (0) \bigcap \mathbb{H}$ is a knot of the type $\lambda$ (bounded or unbounded). Hence, our theorem is proved.
\par\ \par\
$\textbf{\emph{QED}}$
\par\ \par\
Now, given any topological knot class $\lambda$ in $\mathbb{H}$ (bounded or unbounded of finite type), these exists an algebraic knot $K = \{p=0\} \bigcap \mathbb{H}$ which is of that given type $\lambda$. However, as $p \in \mathcal{L}_\mathbb{H}(\mathcal{P})$ (it is onto), there exists some complex polynomial $g$ so that
$\mathcal{L}_\mathbb{H}(g) = p$. Hence, the set of complex tangents to $G = graph(g|_\mathbb{H})$ is precisely
\par\
$\aleph_g = \{\mathcal{L}_\mathbb{H}(g) = 0\} = p^{-1}(0) \bigcap \mathbb{H}$, \par\
the algebraic knot of the arbitrary
(bounded/unbounded) type $\lambda$. We have thus proved the following theorem:
\par\ \par\
\begin{theorem} Every topological type of knot which is either bounded or unbounded of finite count in $\mathbb{H}$  $(\cong \mathbb{R}^3)$ may be assumed precisely as the set of complex tangents to some embedding of $\mathbb{H} \hookrightarrow \mathbb{C}^3$; in fact we can take such as a polynomial embedding and the corresponding knot algebraic.               \end{theorem}
\par\ \par\
We note that our above arguments could not possibly generalize to those unbounded knots of infinite count.
\par\ \par\
$\emph{\textbf{\underline{Note:}}}$ If we were to extend the theorem of Akbulut-King to general links on $S^3$ as we indicated in the previous section, our result could be extended directly to show that (almost) any class of link in $\mathbb{H}$ also admits an algebraic representative. In particular, any link class consisting of all bounded knots except possibly one unbounded knot of finite count. Hence, we would arrive at the fact that every topological type of link in $\mathbb{H}$ satisfying the above condition may arise as the set of complex tangents to some (smooth) embedding $\mathbb{H} \hookrightarrow \mathbb{C}^3$.
\par\ \par\
Further, given any (real) algebraic surface in $\mathcal{A} \subset \mathbb{H}$ given by a real polynomial equation $\{q=0\}$, we change coordinates to $\{z, w, \overline{z}, \overline{w}\}$ coordinates and we can apply Lemma 3 to find that there exists a (polynomial) embedding $\mathbb{H} \hookrightarrow \mathbb{C}^3$ whose complex tangents arise exactly along $\mathcal{A} \subset \mathbb{H}$. Now, with the fact that every topological type of surface, bounded and unbounded, in $\mathbb{R}^3$ (and hence in $\mathbb{H}$) assumes an algebraic representative (see Narasimhan in [8]), we can follow the analogous argument as above for knots to find that every topological type of surface may be assumed
as the set of complex tangents to some (polynomial) embedding of $\mathbb{H} \hookrightarrow \mathbb{C}^3$ (note such a situation is degenerate or rather, not generic).
\par\ \par\

\section*{4. Extension of Results to the 3-sphere}
The tangent CR-operator of $S^3$ is (spanned by) $\mathcal{L}_{S^3}=\mathcal{L}=w\frac{\partial}{\partial\overline{z}}-z\frac{\partial}{\partial\overline{w}}$. It operates on functions which are smooth on $\mathbb{C}^2$.
\par\
Following the same methodology as before, any such map $f$ gives an embedding $S^3 \hookrightarrow \mathbb{C}^3$ via $F = graph(f|_{S^3})$ so that the set of complex tangents of $F(S^3)$ is precisely $\aleph_f = \{\mathcal{L} (f)(z, w) =0\} \bigcap S^3$.
\par\
One may be persuaded to investigate the operator $\mathcal{L}$ and try to directly demonstrate that every topological type of knot may arise as the zero set of some function in the range of $\mathcal{L}$. We do not know how to solve the problem in such generality.
\par\
Here we will attempt to get some results for the 3-sphere by using our above results for the Heisenberg group. Although we do not quite get the strength
of a result as we may have hoped when initially considering the problem, we will show that every topological knot type  in $S^3$ may arise as precisely the set of complex tangents to an embedding into $\mathbb{C}^3$. Unfortunately, we cannot (as of yet) show that we can construct all topological types of knots may arise as complex tangents via means of a $\emph{polynomial}$ embedding, i.e. an embedding given as the graph of a complex polynomial (in two variables), nor are we yet able to show that every knot may arise consisting of totally "non-degenerate" complex tangents.
\par\ \par\
Note that some knots, in particular all torus knots, may arise as complex tangents in such a manner. Let $p,q \in \mathbb{N}$ be relatively prime positive integers. Then the embedding given by the graph of the function $f(z,w)=w^{q-1} \overline{z}-z^{p-1} \overline{w}$ will have its complex tangents precisely when $\mathcal{L} (f) = 0$. But $\mathcal{L} (f) =z^p+w^q$, and so $\aleph_f = \{z^p+w^q=0\} \bigcap S^3$, which is well known to be the standard torus knot
of type $(p,q)$; see Milnor ([7]). Hence, the set of complex tangents to $S^3 \hookrightarrow \mathbb{C}^3$ given by $F=graph(f|_{S^3})$ is precisely a torus knot of type $(p,q)$ in $S^3$.
\par\
One can also similarly
show that any knot that can be given as the intersection of $S^3$ with a complex algebraic (holomorphic) hypersurface passing through the origin may arise as the set of complex tangents to an embedding of $S^3$; this follows as it is easy
to see that all holomorphic functions (with zero constant term) are in the range of $\mathcal{L}$. To demonstrate this, take any such holomorphic function $h(z,w)$ which we can write as: $h(z,w) = z h_1 (z,w) + w h_2 (z,w)$, for some holomorphic functions $h_1, h_2$ (since $h$ has no constant term). Now consider the function $f(z,w) = \overline{z} h_2 (z,w) - \overline{w} h_1 (z,w)$. One readily sees that $\mathcal{L}(f) = h$, and as such the complex tangents to the embedding $F = graph(f|_{S^3})$ will be precisely the zero set of the arbitrary (holomorphic) function $h$ intersected with $S^3$.
\par\ \par\
In the following theorem we demonstrate that there exist embeddings which are complex tangent precisely along a knot of any given type in $S^3$.
\par\
\begin{theorem} For every knot $K \subset S^3$ and integer $n \in \mathbb{N}$ there exists an embedding $Q:S^3 \hookrightarrow \mathbb{C}^3$ of class $\mathcal{C}^n$, so that the set of complex tangents to $Q(S^3)$ is precisely a knot which is equivalent (up to homeomorphism) to $K$; i.e. every topological knot type in $S^3$ may arise precisely as the set of complex tangents to some embedding into $\mathbb{C}^3$. Furthermore, the embedding can be taken to be $\mathcal{C}^\infty$ away from a certain (degenerate) point.    \end{theorem}
\par\ \par\
$\textbf{\emph{\underline{Proof}:-}}$
\par\
Let $K \subset S^3$ be a general knot, and let $\lambda$ be it's topological type.
Assume without loss of generality that $K$ passes through the point $(0,1) \in S^3$ (we may ensure this through a linear rotation). Consider then the knot $\psi(K \setminus \{(0,1)\}) \subset \mathbb{H}$, where $\psi: S^3 \setminus \{(0,1)\} \rightarrow \mathbb{H}$ is the given biholomorphism. Note by our previous arguments and construction before the statement of Theorem 10 in section three that $\psi(K \setminus \{(0,1)\})$ will be a knot of unbounded type $\lambda$ in $\mathbb{H}$.
\par\
Furthermore, by Theorem 12 there exists a knot $\widetilde{K} \subset \mathbb{H}$ so that $\widetilde{K}$ is algebraic and of unbounded knot type $\lambda$ (i.e., topologically equivalent to $\psi(K \setminus \{(0,1)\})$). In particular, let $g:\mathbb{C}^2 \rightarrow \mathbb{C}$ be the complex polynomial so that
$\widetilde{K} = \{g=0\} \bigcap \mathbb{H}$. Then, as $\mathcal{L}_\mathbb{H}$ is an onto operator on polynomials (by Lemma 9),
there exists a polynomial $f:\mathbb{C}^2 \rightarrow \mathbb{C}$ so that: $\mathcal{L}_\mathbb{H}(f)=g$.
\par\
Hence, the set of complex tangents to $F=graph(f|_\mathbb{H}):\mathbb{H} \hookrightarrow \mathbb{C}^3$ is precisely $\widetilde{K}$.
\par\ \par\
Let $\mathcal{M}_f = F(\mathbb{H}) =\{(x,f(x)) | x \in \mathbb{H} \} \subset \mathbb{C}^3$. Further, let $\mathcal{S}_{f \circ \psi}$ be the space given as the graph of:
$f \circ \psi: S^3 \setminus \{(0,1)\} \rightarrow \mathbb{C}$.
\par\
Note then: $\mathcal{S}_{f \circ \psi} = \{ (s, f( \psi(s)) | s \in S^3\ \setminus \{(0,1)\} \} \subset \mathbb{C}^3$. But then clearly $\mathcal{M}_f$ and $\mathcal{S}_{f \circ \psi}$ are biholomorphic, in particular it will be given by: $\varphi \times Id:\mathcal{M}_f \rightarrow \mathcal{S}_{f \circ \psi}$ where $\varphi=\psi^{-1}$ acts on the manifold coordinates and identity on the (last) complex coordinate. Clearly this map will be holomorphic and diffeomorphic (as $\varphi$ and $Id$ are), and so is its inverse.
\par\
Therefore, as biholomorphisms preserve complex tangents, $\mathcal{S}_{f \circ \psi}$ and $\mathcal{M}_f$ will have equivalent sets of complex tangents. In fact, it is clear that the set of complex tangents to $\mathcal{S}_{f \circ \psi}$ will be precisely $\varphi(\widetilde{K})$, since its the biholomorphic image of the set of complex tangents to $\mathcal{M}_f$, which form $\widetilde{K}$.
\par\ \par\
Hence, we have an embedding $S^3 \setminus \{(0,1)\} \hookrightarrow \mathbb{C}^3$ whose complex tangents will form a knot of type $\lambda$ in $S^3$ if we were to add the point at infinity (0,1). However, this embedding is given precisely as the graph of the rational function $f \circ \psi$ over $S^3 \setminus \{(0,1)\}$. As $ \psi(z,w) = (\frac{iz}{1-w},i\frac{1+w}{1-w})$ and $f$ is a complex polynomial, we see that $f \circ \psi$ is a smooth function away from $w=1$, i.e. the point $(0,1) \in S^3$. In fact, as we approach $(0, 1)$, the function becomes unbounded. Hence, we will
not be able to directly extend this function to the entire sphere (so as to garner an embedding of the entire sphere).  \par\
Suppose $f$ is a complex polynomial of degree $n$. Then: \par\
$(f \circ \psi) (z, w) = f(\frac{iz}{1-w},i\frac{1+w}{1-w})$, which we may then write: $(f \circ \psi) (z, w) = \frac{1}{(1-w)^n (1- \overline{w})^n)} p(z,w)$, where $p$ is a complex polynomial. This can be done as in each term of $(f \circ \psi) (z, w)$ the denomonator will consist precisely of a term: $\frac{1}{(1-w)^k (1- \overline{w})^l}$, where  $k+l \leq n$ as $f$ is a polynomial of degree n. Hence, multiplying $f$ by $(1-w)^{2n+r}$
(for any $r \in \mathbb{N}$) gives us:
\par\
${q_f}^r=(1-w)^{2n+r} f(z,w) = \frac{(1-w)^{n+r}}{(1-\overline{w})^n} p(z,w)$,
\par\
where $p$ again is a complex polynomial.
\par\ \par\
Now, it is well known (and easy to see using simple calculus) that the rational function in one complex variable: $h(\zeta)=\frac{\zeta^k}{\overline{\zeta}^l}$
will be everywhere smooth except at the point $\zeta=0$. Furthermore, if $k=l+r$ where $k, l, r \in \mathbb{N}$ (positive integers), then $f$ will be
of class $\mathcal{C}^{r-1}$ at $\zeta=0$, i.e. $f$ will be $(r-1)$-times continuously differentiable at $0$, and smooth (of class $\mathcal{C}^\infty$)
everywhere else.
\par\
Now in our situation, ${q_f}^r = \frac{(1-w)^{n+r}}{(1-\overline{w})^n} p(z,w)$ and $p$ is a polynomial. By our above remark, (considering $\zeta=1-w$) we see that the rational term $\frac{(1-w)^{n+r}}{(1-\overline{w})^n}$ will be continuous at $(0,1)$ if $r \geq 1$ and its value there will be $0$, and further for larger $r$ it will be $(r-1)$-times continuously differentiable at $(0,1)$.
\par\ \par\
Take $r \geq 2$ (so that ${q_f}^r$ is continuously differentiable). Let $\mathcal{Q}_r = Q_r(S^3) \subset \mathbb{C}^3$,
where $Q_r = graph({q_f}^r |_{S^3})$. Then $\mathcal{Q}_r$ will be homeomorphic to $S^3$, and whose tangent spaces are well-defined and vary continuously.
\par\
Note that we have a map: $\Upsilon:\mathcal{S}_{f \circ \psi} \rightarrow \mathcal{Q}_r \setminus \{(0,1,0)\}$ given by:
$(s, f(\psi(s)) \rightarrow (s, (1+w)^{2n+r} f(\psi(s))$. This map may be given precisely by: $Id \times \tau$, where $\tau(\zeta)=(1+w)^{2n+r} \zeta$. Note this map is one-to-one (it is identity of first coordinate, and the manifold is a graph) and onto by construction,
and furthermore it is a holomorphism as it is merely given by multiplication with a holomorphic function in the second coordinate. It is also easily checked that $\Upsilon^{-1}$ is a holomorphism, as it is given by multiplication by $\frac{1}{(1+w)^{2n+r}}$ in the second coordinate which is also holomorphic (away from (0,1)).
\par\ \par\
Hence, $\Upsilon$ is a biholomorphism and so $\mathcal{Q}_r \setminus \{(0,1,0)\}$ will have the same complex tangents as $\mathcal{S}_{f \circ \psi}$ over
$S^3 \setminus \{(0,1,0)\}$. The set $\varphi(\widetilde{K})$ constitutes the complex tangents and forms a "punctured" knot of type $\lambda$ (that is, it is missing $(0,1)$ to be a closed knot).
\par\
Now, it remains to resolve the situation at the point $(0,1) \in S^3$, or equivalently at $(0,1,0) \in \mathcal{Q}_r$. Note that when we apply
 $\mathcal{L}$ (the CR-operator of $S^3$) to ${q_f}^r$, the term $(1-w)^{n+r}$ factors out:
\par\ \par\
 $\mathcal{L}({q_f}^r) (z,w) = \mathcal{L}(\frac{(1-w)^{n+r}}{(1-\overline{w})^n} p(z,w)) = \frac{(1-w)^{n+r}}{(1-\overline{w})^n} \mathcal{L}(p) (z,w) - n \frac{(1-w)^{n+r}}{(1-\overline{w})^{n+1}} z p(z,w) $.
\par\
As $r \geq 2$, both of the terms: $\frac{(1-w)^{n+r}}{(1-\overline{w})^n}, \frac{(1-w)^{n+r}}{(1-\overline{w})^{n+1}}$ are continuous at $(0,1)$ and take a value of $0$ at such points. Hence $\mathcal{L}({q_f}^r) (0,1) = 0$ and $(0,1) \in S^3$ is a complex tangent of the embedding $\mathcal{Q}_r$.
But we know that away from $(0,1)$ the complex tangents of $\mathcal{Q}_r$ form a knot of type $\lambda$ punctured at $(0,1)$
\par\ \par\
Hence, we immediately find: $\aleph_{\mathcal{Q}_r} = \varphi(\widetilde{K}) \bigcup \{(0,1)\} = C_r$,
 which is precisely a knot of type $\lambda$, smooth everywhere but $(0,1)$, but for any given $r \geq 2$ the knot will $(r-1)$-times differentiable at the point $(0,1)$. As $\lambda$ was an arbitrary equivalence class of knots, we have proved the claim of the theorem.
\par\ \par\
$\textbf{\emph{QED}}$
\par\ \par\
The constructed knots will be $\mathcal{C}^\infty$ smooth everywhere but at the point $(0,1)$, where they can be taken to be continuously differentiable $n$-times (for any given $n$). This point will be degenerate, as all of the derivatives (that are defined) at the point will be zero. One can also easily check that the Bishop invariant must also be degenerate by construction; see our forthcoming paper [3]. As a result, an arbitrary perturbation of the embedding will give rise to a more complicated set of complex tangents (in general). In particular, we may expect to get a knot which is  equivalent to the knot $K$ away from $(0,1)$, but near $(0,1)$ the topology may be drastically altered and we may find new knot components in a neighborhood of $(0,1)$. A potential direction of future work would be to create a "smoothing" argument that will allow us to generate a smooth embedding without adding any more knot components or changing the knot topology of $K$.
\par\ \par\
We note that we can apply the above construction to a knot which doesn't pass through the point $(0,1)$.  This will allow us to construct embeddings of the 3-sphere whose complex tangents will form precisely a smooth ($\mathcal{C}^\infty$) knot of the same type as the given knot and will be smooth away from the point $(0,1)$. However, in doing so we will generate a single complex tangent point at $(0,1)$, which will be again degenerate, and the embedding will be only of class $\mathcal{C}^n$ near the point. If one were to perturb such an embedding (by a small perturbation), we would in general get the knot $K$ union some other knot (or link) in the neighborhood of the point. We are as of yet unable to "remove" this complex tangent point at infinity to get the desired result for smooth embeddings.
\par\ \par\
We will now extend our theorem directly to simple chain links, i.e. links whose components are knots  $\{K_1,...,K_n\} \subset S^3$ arranged in such a fashion that $K_j$ links with $K_{j-1}$ and $K_{j+1}$ with linking number one for $2 \leq j \leq n-1$, and $K_1$ links with $K_2$ and $K_n$ links with $K_{n-1}$, both also of linking number one. We give the proof of our desired result (Theorem 14, see below) by the following construction.
\par\ \par\
Let $\{K_1,...,K_n\} \subset S^3$ be of respective types: $\{\lambda_1,...,\lambda_n\}$ (some possibly equal) and assume without loss of generality that $K_1,...,K_n$ are disjoint and that $(0,1) \in K_1$.
\par\
Map these knots using $\psi$ to the Heisenberg group, to get $n$ knots $\psi(K_1),...,\psi(K_n)$ of respective types $\{\lambda_1,...,\lambda_n\}$. They are bounded knots except $\psi(K_1)$ is of unbounded type. Again, by our previous results, there exists algebraic knots $\widetilde{K}_i\subset \mathbb{H}$ so that $\widetilde{K}_1$ is of unbounded type $\lambda_1$ and $\widetilde{K}_j$ is of bounded type $\lambda_j$, for each $2 \leq j \leq n$.
\par\ \par\
Now we will need to "move" these knots in an algebraic way so that they are simple chain links of the analogous configuration to that described above in $S^3$. Let $c \in \mathbb{H}$ and let $\Pi:\mathbb{H} \rightarrow \mathbb{R}^3$ by the  diffeomorphism we previously mentioned, with $\Pi(z, w) = (Re(z), Im(z), Re(w)) = (x, y, z) \in \mathbb{R}^3$.
Let $\tau_a:\mathbb{R}^3 \rightarrow \mathbb{R}^3$ be translation by a vector $a$, i.e. $\tau_a (x) = x+a$. Note that $\tau_a$ is a diffeomorphism.
Now consider the diffeomorphism: $\widetilde{\tau}_c = \Pi^{-1} \circ \tau_{\Pi(c)} \circ \Pi: \mathbb{H} \rightarrow \mathbb{H}$.
Note that this is a diffeomorphism with $\widetilde{\tau}_c (0,0) = c$, and $\widetilde{\tau}_c$ will obviously preserve algebraic sets.
Hence, $\widetilde{\tau} (\widetilde{K})$ will be algebraic knot of the same type $\lambda$ (boundedness also preserved) but translated from $0$ to $c$.
\par\ \par\
Therefore we may translate our knots without sacrificing their type and keeping them algebraic. Further, note that given any linking formation of knots in $S^3$ there is an analogous arrangement in $\mathbb{H}$ (by formulation).
\par\
Note also that we may rotate knots in $\mathbb{H}$ in a linear manner, namely by multiplying by $(3 \times 3)$-matrices in $\mathbb{R}^3$ and
(pre-)composing with $\Pi, \Pi^{-1}$. Note that such operations must preserve algebraic sets (they are all "polynomial transformations").
\par\
Now given the implicit simple chain linking arrangement in $\mathbb{H}$, it is clear with some geometric intuition that
we may translate and rotate the knots $\widetilde{K}_j$ so they link to form the desired chain, as they link in a chain with linking number one and they are 1-manifolds in 3-space.
\par\ \par\
Assume then without loss of generality that the $\widetilde{K}_j$ link in the appropriate (simple chain) manner and that their types are of  $\{\lambda_1,...,\lambda_n\}$ respectively, with only $\widetilde{K}_1$ being unbounded. Let $g_i$ be complex polynomials whose zero sets are precisely $\widetilde{K}_i$, respectively. Note that as the knots are disjoint, if we take $g_* = g_1 g_2 ... g_n$, the zero set of $g_*$ will be our desired link (when restricted to $\mathbb{H}$).
\par\
Let $f$ be a polynomial so that $\mathcal{L}_\mathbb{H}(f) = g_*$. Then, we can follow our argument in the proof of Theorem 13 with this new function $f$ to find an embedding of $S^3$ whose complex tangents are precisely a simple chain link biholomorphic to the set $\{g_* = 0\} \subset \mathbb{H}$ (union the north pole). By the tautological relationship (by definition) of types of links in $S^3$ with types of links in $\mathbb{H}$, we obtain the result below:
\begin{theorem}: For any simple chain link of knots $L \subset S^3$, there exists an embedding $S^3 \hookrightarrow \mathbb{C}^3$ whose complex tangents form a link of knots that is topologically equivalent to $L$.     \end{theorem}
\par\ \par\
$\textbf{\emph{\underline{Remark:}}}$ Our above construction (via translation and rotation) would not apply for general links. However if one were to generalize the work of Akbulut-King to general links on $S^3$ as we indicated in section three, we could then generalize our proof to "all links are algebraic" in $\mathbb{H}$. We would then immediately obtain the desired result on $S^3$ for links; namely our proof for Theorem 13 will hold for links and we will arrive at the fact that given any topological class of link in $S^3$ there is a link of that class that is assumed exactly as the set of complex tangents to some $\mathcal{C}^k$ embedding (with still one degenerate point).
\par\ \par\
We also note that we can extend our theorem above to show that any algebraic surface $\mathcal{A} \subset S^3$ may arise as the set of complex tangents to some embedding $S^3 \hookrightarrow \mathbb{C}^3$, using the corresponding result for the Heisenberg group (see our note at the end of section three). Assuming any topological type of surface may also arise in $S^3$ as an algebraic set (and similarly in $\mathbb{H}$), we may follow our argument above for knots in exact analogy to establish that for every  surface $D \subset S^3$, there exists an embedding $S^3 \hookrightarrow \mathbb{C}^3$ whose complex tangents form a surface which is topologically equivalent to $D$.
\par\
Furthermore, taking products of polynomials (as above), we may achieve any simple union of such topological types (surfaces or knots) as the set of complex tangents to some embedding of $S^3$. We demonstrate some examples of interesting (and singular) possible scenarios in our paper [3]. We refer the reader to this paper for the computation of some examples and an analysis of our degenerate point.
\par\ \par\
We may also generalize some of our work above and get interesting results for (higher) odd-dimensional spheres. We will exhibit this in a future paper.
\par\ \par\ \par\ \par\
\par\ \par\ \par\

\end{document}